\documentclass[reqno]{amsart}
\usepackage{amsmath,amscd,amsfonts,amssymb}
\usepackage{mathrsfs,dsfont}
\usepackage{hyperref}

\def\R{\mathbb{R}}

\def\Z{\mathbb{Z}}
\def\T{\mathbb{T}}
\def\1{\mathds{1}}

\newcommand\eqdef{\overset{\text{def}}{=}}

\newcommand\gab{\Z \alpha + \Z \beta}

\renewcommand\leq{\leqslant}

\newcommand{\ft}[1]{\widehat #1}
\newcommand\dotprod[2]{\langle #1 , #2 \rangle}
\newcommand\mes{\operatorname{mes}}
\newcommand\const{\operatorname{const}}

\theoremstyle{plain}
\newtheorem{theorem}{Theorem}

\newtheorem{lemma}{Lemma}

\theoremstyle{definition}

\newenvironment{enumerate-math}
{\begin{enumerate}
\addtolength{\itemsep}{5pt}
}
{\end{enumerate}}

\newenvironment{enumerate-math-abc}
{\begin{enumerate}
\addtolength{\itemsep}{5pt}
}
{\end{enumerate}}

\begin{document}

\title{Riesz bases of exponentials on multiband spectra}

\author{Nir Lev}
\address{Department of Mathematics, Weizmann Institute of Science, Rehovot 76100, Israel.}
\email{nir.lev@weizmann.ac.il}

\subjclass[2000]{42C15, 94A12}

\keywords{Riesz bases, Multiband signals, Quasicrystals}

\begin{abstract}
Let $S$ be the union of finitely many disjoint intervals on $\R$.
Suppose that there are two real numbers $\alpha, \beta$ such that
the length of each interval belongs to $\gab$. We use quasicrystals
to construct a discrete set $\Lambda \subset \R$ such that the
system of exponentials $\{\exp 2 \pi i \lambda x, \, \lambda \in
\Lambda\}$ is a Riesz basis in the space $L^2(S)$.
\end{abstract}

\maketitle


\section{Introduction}

\subsection{}
Let $S$ be the union of a finite number of bounded intervals on $\R$.
We denote by $PW_S$ (after Paley and Wiener) the space of all functions
$f \in L^2(\R)$ whose Fourier transform
\[
\ft{f}(x) = \int_{\R} f(t) \, e^{2 \pi i x t} \, dt
\]
vanishes almost everywhere outside of $S$. A discrete set $\Lambda \subset \R$ is
called a \emph{complete interpolation set} for $PW_S$ if the restriction
operator $f \mapsto f|_\Lambda$ is a bounded and invertible one from
$PW_S$ onto $\ell^2(\Lambda)$. In the context of communication theory
this means that $\Lambda$ provides a ``stable and non-redundant''
sampling of signals with spectrum in $S$.

It is well-known that the complete interpolation property of $\Lambda$
is equivalent to the Riesz basis property of the corresponding exponential
system
\[
E(\Lambda) = \{ \exp 2 \pi i \lambda x, \; \lambda \in \Lambda\}
\]
in the space $L^2(S)$.

If $S$ is a single interval, then a complete description of the Riesz bases
$E(\Lambda)$ in $L^2(S)$ was given by B.\ S.\ Pavlov (1979). Much less is
known, however, in the case when $S$ is the union of more than one interval.
In fact, it is unknown in general whether an exponential Riesz basis in 
$L^2(S)$ exists at all. This existence has been established in the following
special cases:
\begin{enumerate-math}
\item
$S$ is a finite union of disjoint intervals with commensurable lengths
\cite{bezuglaya-katsnelson, lyubarskii-seip}.
\item
$S$ is the union of two general intervals \cite{seip}.
\end{enumerate-math}
For other results in the subject we refer to the survey paper
\cite{lyubarskii-seip}.

In this note we extend the two results above and prove:
\begin{theorem}
\label{thm:main-1}
Let $S$ be the union of finitely many disjoint intervals on $\R$.
Suppose that there are two real numbers $\alpha, \beta$ such that
the length of each interval belongs to $\gab$. Then there is $\Lambda
\subset \R$ such that $E(\Lambda)$ is a Riesz basis in $L^2(S)$.
\end{theorem}
Here $\gab$ denotes the set of real numbers of the form
$n \alpha + m \beta$ $(n,m \in \Z)$. The above mentioned results are
thus obtained as special cases of Theorem \ref{thm:main-1}.

\subsection{}
In fact, we will prove the following more general result.

\begin{theorem}
\label{thm:main-2}
Suppose that the indicator function of a set $S \subset \R$ can be
expressed as a linear combination of indicator functions of intervals
$I_1, \dots, I_N$ whose lengths belong to $\gab$, that is,
\begin{equation}
\label{eq:linear-combination}
\1_S(x) = \sum_{j=1}^{N} c_j \, \1_{I_j}(x), \quad
|I_j| \in \gab \quad (1 \leq j \leq N),
\end{equation}
where $\alpha, \beta \in \R$. Then there is $\Lambda \subset \R$
such that $E(\Lambda)$ is a Riesz basis in $L^2(S)$.
\end{theorem}

Sets with the structure \eqref{eq:linear-combination} form a wider
class than unions of disjoint intervals with lengths in $\gab$.
For example, one may take an interval with ``holes'' obtained by 
the removal of disjoint sub-intervals, where the interval and its
sub-intervals have their lengths in $\gab$.

A result of similar type in the periodic setting was obtained in \cite{kozma-lev}.


\section{Quasicrystals. Duality.}

Our approach is inspired by the papers \cite{matei-meyer-quasicrystals,
matei-meyer-simple} due to Matei and Meyer, who introduced the usage
of so-called `simple quasicrystals' in order to construct ``universal''
sets of sampling or interpolation for $PW_S$ spaces. Here we will use
simple quasicrystals to construct complete interpolation sets for spectra
$S$ with the structure \eqref{eq:linear-combination}.

Following \cite{matei-meyer-quasicrystals, matei-meyer-simple} we let
$\Gamma$ be a lattice in $\R^2$. Consider the projections $p_1(x,y) = x$
and $p_2(x,y)=y$, and assume that the restrictions of $p_1$ and $p_2$ to
$\Gamma$ are injective. Let $\Gamma^*$ be the dual lattice, consisting
of all vectors $\gamma^* \in \R^2$ such that $\dotprod{\gamma}{\gamma^*}
\in \Z$, $\gamma \in \Gamma$.

Let $S$ be the union of disjoint semi-closed intervals,
\begin{equation}
\label{eq:def-s}
S = \bigcup_{j=1}^{\nu} [a_j, b_j), \quad a_1 < b_1 < \cdots < a_\nu < b_\nu,
\end{equation}
and
\begin{equation}
\label{eq:def-i}
I = [a, b)
\end{equation}
be a single semi-closed interval. Define
\[
\begin{aligned}
\Lambda(\Gamma, I) &= \{ p_1(\gamma) : \gamma \in \Gamma, \; p_2(\gamma) \in I\},\\[4pt]
\Lambda^*(\Gamma, S) &= \{ p_2(\gamma^*) : \gamma^* \in \Gamma^*, \; p_1(\gamma^*) \in S\}.
\end{aligned}
\]

In some sense, the quasicrystals $\Lambda(\Gamma, I)$ and $\Lambda^*(\Gamma, S)$
are dual to each other. This duality was observed and used by Matei and Meyer
in connection with sampling and interpolation. In the present context of
exponential Riesz bases the duality can be formulated as follows.

\begin{lemma}
\label{lemma:duality}
The following two properties are equivalent:
\begin{enumerate-math}
\item
$E(\Lambda(\Gamma, I))$ is a Riesz basis in $L^2(S)$;
\item
$E(\Lambda^*(\Gamma, S))$ is a Riesz basis in $L^2(I)$.
\end{enumerate-math}
\end{lemma}

The proof of Lemma \ref{lemma:duality} is along similar lines as in
the paper \cite{matei-meyer-simple} (Sections 6--7), but in our case
there is an additional point concerned with the requirement that the
intervals in \eqref{eq:def-s} and \eqref{eq:def-i} should be semi-closed. 
Indeed, the significance of this point is clarified once keeping
in mind that the Riesz basis property of the exponential systems
$E(\Lambda(\Gamma, I))$ and $E(\Lambda^*(\Gamma, S))$ is not
preserved upon either the addition or removal of any element.

For a proof of the duality lemma in the periodic setting
see \cite[Section 2]{kozma-lev}.


\section{Proof of Theorem \ref{thm:main-2}}

There is no loss of generality in assuming that the numbers $\alpha, \beta$ are
linearly independent over the rationals. Moreover, by rescaling we may restrict
ourselves to the case when
\begin{equation}
\label{eq:irrational}
\text{$\alpha$ is an irrational number, and $\beta = 1$.}
\end{equation}

\subsection{}
Define a lattice
\[
\Gamma = \{(n(1 + \alpha) - m, \, m - n \alpha) : n,m \in \Z\},
\]
and let $I = [0, \mes S)$ be an interval whose length coincides with the Lebesgue
measure of $S$. We will prove that the exponential system $E(\Lambda(\Gamma, I))$
is a Riesz basis in $L^2(S)$. According to Lemma \ref{lemma:duality} it will be
sufficient to show that the system $E(\Lambda^*(\Gamma, S))$ is a Riesz basis in $L^2(I)$.

It is easy to check that the set $\Lambda^*(\Gamma, S)$ may be partitioned
as follows,
\begin{equation}
\label{eq:partition}
\Lambda^*(\Gamma, S) = \bigcup_{n \in \Z} \Lambda_n, \quad
\Lambda_n = (S \cap (n \alpha + \Z)) + n
\end{equation}
(where some of the sets $\Lambda_n$ may be empty). Let $\{s_n\}$
be a sequence of integers such that $s_n - s_{n-1} = \# \Lambda_n$, and choose
an enumeration $\{\lambda_j, \, j \in \Z\}$ of the set $\Lambda^*(\Gamma, S)$
such that
\[
\Lambda_n = \{ \lambda_j \, : \, s_{n-1} \leq j < s_n\} \quad (n \in \Z).
\]

In order to prove that $E(\Lambda^*(\Gamma, S))$ is a Riesz basis in $L^2(I)$
it will be sufficient, by a theorem of Avdonin \cite{avdonin}, to check that
the following three conditions hold:
\begin{enumerate-math-abc}
\item
\label{item:separation}
$\{\lambda_j\}$ is a separated sequence, $\inf_{j \neq k} |\lambda_j - \lambda_k| > 0$;
\item
\label{item:boundedness}
$\sup_{j} |\delta_j| < \infty$, where $\delta_j = \lambda_j - j/\mes S$;
\item
\label{item:cancellation}
There is a constant $c$ and a positive integer $N$ such that
\begin{equation}
\label{eq:avdonin-condition}
\sup_{a \in \Z} \; \Big| \frac1{N} \sum_{j=a+1}^{a+N} \delta_{j} \; - \; c\Big| < \frac1{4 \, \mes S} \, .
\end{equation}
\end{enumerate-math-abc}

Condition \ref{item:separation} can be easily verified directly from the 
definition of $\Lambda^*(\Gamma, S)$.

\subsection{}
We will next show that \ref{item:boundedness} holds. Consider the bounded,
piecewise constant function
\[
\phi(x) = \sum_{k \in \Z} \1_S(x - k).
\]
This function is $1$-periodic and hence may be viewed as a function on the
circle group $\T = \R / \Z$. The assumptions \eqref{eq:linear-combination},
\eqref{eq:def-s} and \eqref{eq:irrational} now provide the following
representation for $\phi$.

\begin{lemma}
\label{lemma:oren}
There is a bounded, piecewise linear function $\psi: \T \to \R$ such that
\[
\phi(x) = \mes S + \psi(x) - \psi(x - \alpha), \quad x \in \T.
\]
\end{lemma}

For a proof of Lemma \ref{lemma:oren} see \cite[Lemma 3.2]{kozma-lev}.
Now observe that by \eqref{eq:partition} we have
\begin{equation}
\label{eq:sn-phi}
s_n - s_{n-1} = \# \Lambda_n = \phi(n \alpha).
\end{equation}
It follows from \eqref{eq:sn-phi} and Lemma \ref{lemma:oren} that
\begin{equation}
\label{eq:sn-formula}
s_n = n \, \mes S + \psi(n \alpha) + \const.
\end{equation}
Given $j$ there is $n = n(j)$ such that $\lambda_j \in \Lambda_n$, or equivalently,
such that $s_{n-1} \leq j < s_n$. Then
\[
\delta_j = \lambda_j - \frac{j}{\mes S}
= (\lambda_j - n) + \Big(n - \frac{s_n}{\mes S}\Big) + \frac{s_n - j}{\mes S}.
\]
It thus follows from \eqref{eq:partition}, \eqref{eq:sn-phi} and \eqref{eq:sn-formula}
that $\sup_j |\delta_j| < \infty$, which confirms condition \ref{item:boundedness} above.

\subsection{}
It remains to establish \ref{item:cancellation}. In order to show that
\eqref{eq:avdonin-condition} holds we will first obtain a simple expression
for the sum $\sum \delta_j$ where $j$ goes through the interval
$s_{n-1} \leq j < s_n$. Indeed,
\[
\sum_{j=s_{n-1}}^{s_n - 1} \delta_j
= \sum_{j=s_{n-1}}^{s_n - 1} (\lambda_j - n)
- \sum_{j=s_{n-1}}^{s_n - 1} \Big(\frac{j}{\mes S} - n \Big)
\eqdef S_1(n) - S_2(n).
\]
We evaluate each one of the sums $S_1(n),S_2(n)$ separately. First we observe
that by \eqref{eq:partition},
\[
S_1(n) = \sum_{k \in \Z} (n \alpha - k) \, \1_S(n \alpha - k) \eqdef \tau_1(n \alpha).
\]
Secondly, by a direct calculation and using \eqref{eq:sn-phi} and 
\eqref{eq:sn-formula} we find that
\[
S_2(n) = (s_n - s_{n-1}) \Big( \frac{s_{n-1} + s_n - 1}{2 \mes S} - n \Big)
= \frac{\phi(n \alpha) (\psi(n \alpha) - \frac1{2} \phi(n \alpha) + \const)}{\mes S}
\eqdef \tau_2(n \alpha).
\]
We conclude that for an appropriately defined function $\tau : \T \to \R$
(bounded and piecewise continuous) we have
\begin{equation}
\label{eq:sum-simple}
\sum_{j=s_{n-1}}^{s_n - 1} \delta_j = \tau(n \alpha), \quad n \in \Z.
\end{equation}

\subsection{}
Now we can finish the proof of \ref{item:cancellation} above. Given $a \in \Z$
and a positive (large) integer $N$ there are $n = n(a)$ and $r = r(a, N)$ such that
\[
s_{n-1} \leq a < s_n, \quad s_{n+r-1} \leq a + N < s_{n+r}.
\]
Since the sequence $\{\delta_j\}$ is bounded and due to \eqref{eq:sum-simple} we have
\[
\sum_{j=a+1}^{a+N} \delta_{j} = \sum_{j=s_{n-1}}^{s_{n+r-1} - 1} \delta_{j} + O(1)
= \sum_{k=n}^{n+r-1} \tau(k \alpha) + O(1).
\]
The points $\{n \alpha\}$ are well-distributed on the circle $\T$ (since $\alpha$
is irrational) and hence
\[
\sum_{k=n}^{n+r-1} \tau(k \alpha) = r \int_{\T} \tau(x) \, dx + o(r), \quad r \to \infty,
\]
uniformly with respect to $n$. Since \eqref{eq:sn-phi} and \eqref{eq:sn-formula}
imply that $N = r \, \mes S + O(1)$, we get
\[
\frac1{N} \sum_{j=a+1}^{a+N} \delta_{j} = \frac1{\mes S} \int_{\T} \tau(x) \, dx + o(1),
\quad N \to \infty,
\]
with the $o(1)$ uniform with respect to $a$. This implies \ref{item:cancellation},
and so Theorem \ref{thm:main-2} is proved.


\section{Remarks}

We have constructed Riesz bases of exponentials with real frequencies for
multiband spectra subject to the diophantine condition \eqref{eq:linear-combination}.
We would like though to comment on one result which is not covered by our theorems
above. In the paper \cite{seip} existence of such Riesz bases was proved under
certain non-discrete conditions on the lengths of the gaps between the intervals.
The restrictions obtained are rather severe; however, the result indicates that
diophantine restrictions are not necessarily natural ones in the problem.

We also refer the reader to the paper \cite{lyubarskii-spitkovsky} where
the authors construct, for any finite union of intervals, a Riesz basis
of exponentials with complex frequencies lying in a horizontal strip along
the real axis.

\subsection*{Acknowledgement}
We thank Kristian Seip for reading an earlier version of this note.


\end{document}